\documentclass[english]{ccdconf}
\usepackage{graphicx }
\usepackage{amssymb}
\usepackage{mathrsfs}
\usepackage{enumerate}
\usepackage{multirow}
\usepackage{booktabs}
\usepackage[colorlinks,
            linkcolor=blue,
            anchorcolor=blue,
            citecolor=blue]{hyperref}
\usepackage[numbers,sort&compress]{natbib}
\begin{document}

\title{Analytical calculation of the inverse nabla Laplace transform}
\author{Yiheng Wei\aref{ustc}, YangQuan Chen\aref{ucm}, Yuquan Chen\aref{ustc},
      Yong Wang\aref{ustc}$^\ast$}

\affiliation[ustc]{University of Science and Technology of China, Hefei, Anhui 230026, China
        \email{yongwang@ustc.edu.cn}}
\affiliation[ucm]{University of California, Merced, Merced, California 95343, USA}
\maketitle

\begin{abstract}
The inversion of nabla Laplace transform, corresponding to a causal sequence, is considered. Two classical methods, i.e., residual calculation method and partial fraction expansion method are developed to perform the inverse nabla Laplace transform. For the first method, two alternative formulae are proposed when adopting the poles inside or outside of the contour, respectively. For the second method, a table on the transform pairs of those popular functions is carefully established. Besides illustrating the effectiveness of the developed methods with two illustrative examples, the applicability are further discussed in the fractional order case.
\end{abstract}

\keywords{Analytic inversion, Nabla Laplace transform, Residual calculation method, Partial fraction expansion method, Discrete fractional calculus}

\footnotetext{The work described in this paper was supported by the National Natural Science Foundation of China (61601431,
61573332), the Anhui Provincial Natural Science Foundation (1708085QF141) and the funds of China Scholarship Council
(201706340089, 201806345002).}

\section{INTRODUCTION}\label{Section 1}
Nabla Laplace transform is a powerful tool to convert a discrete time signal, which is a sequence of real or complex numbers, into a complex frequency domain representation \cite{Hein:2011PAMJ}. Similar to the classic $Z$ transform, it can be considered as a discrete time generalization of the Laplace transform \cite{Jarad:2012ADE,Ortigueira:2013IFAC,Ortigueira:2016CNSNS}. It is crucial for a wide range of applications, whenever signal sampling or discrete treatment or backward difference are involved. Especially for the fractional order case, nabla Laplace transform technique is much more convenient and effective than $Z$ transform \cite{Abdeljawad:2012ADE,Jonnalagadda:2016IJNS,Abdeljawad:2017RMP}.

The basic idea now known as the nabla Laplace transform was inspired from Bohner and Peterson \cite{Bohner:2001Book}, and it was formally introduced and investigated in 2009 by At{\i}c{\i} and Eloe as a way to treat fractional finite difference equation \cite{Atici:2009EJQTDE}. It gives a tractable way to solve linear time-invariant fractional backward difference equations whose sampling period $h=1$. It was also dubbed ``the $N$-transform'' for short. Afterwards, the sampling period $h$ was generalized as a positive real number, and then a sampling based nabla Laplace transform was developed \cite{Cheng:2011Book}.

For the constructed nabla Laplace transform, some important properties were presented, e.g. linearity, shifting in the time domain, and convolution theorem \cite{Mohan:2014CMS}. Afterwards, several novel properties were derived and then applied in the nabla fractional calculus (see Chapter 3 of the famous monograph \cite{Goodrich:2015Book}). In \cite{Wei:2018ArXiva}, we established the initial/final value theorem, the stable criterion and then applied such properties to analyze the monotonicity and overshoot properties of the zero input system response. To further explore the properties of such a handful tool, a comprehensive survey was made by us on the existing results, and 14 innovative properties were proposed subsequently \cite{Wei:2019FDTA}. By using this tool, six kinds of infinite dimensional frequency distributed models were equivalent derived to describe a nabla fractional order system \cite{Wei:2019CNSNS}. To obtain the time-domain sequence $f(k)$ of a given function $F(s)$ in frequency domain, a rational approximation approach was proposed in \cite{Wei:2019AJC}. Though this method is effective, it always leads to the approximation error. To get the exact value, a method generalized from the initial/final value theorem was derived in \cite{Wei:2018ArXiva}. However, this method need calculate the limit for each $k$ in $f(k)$. To obtain the explicit expression of $f(k)$, we developed the inverse nabla Laplace transform in the form of a contour integral \cite{Wei:2019AJC}. To be honest, it is generally difficult or even impossible to calculate such a contour integral. Additionally, it is well known that the residual calculation method and the partial fraction expansion method perform well in solving the similar problems in inverse Laplace transform and inverse $Z$ transform. Motivated by these, we will develop the two analytical method for inverse nabla transform Laplace transform.

The outline of the rest paper is as follows. In Section \ref{Section 2}, the basic definition on nabla Laplace transform is review. In Section \ref{Section 3}, two methods are designed and discussed to evaluate the contour integral effectively. In Section \ref{Section 4}, two typical examples are provided to show the feasibility and effectiveness of the developed methods. Finally, some conclusions are drawn in Section \ref{Section 5}.

\section{PRELIMINARIES}\label{Section 2}
In this section, some basic definitions and concepts for nabla Laplace transform are provided. Afterwards, the objective of this work is restated.

The nabla Laplace transform of a sequence $f: \mathbb{N}_{a+1}\to \mathbb{R}$ is defined by \cite{Atici:2009EJQTDE}
\begin{equation}\label{Eq1}
{\textstyle
{\mathscr N}_a\left\{ {f\left( k \right)} \right\} \triangleq \sum\nolimits_{k = 1}^{+\infty}  {{{\left( {1 - s} \right)}^{k - 1}}f\left( k+a \right)},}
\end{equation}
where $s\in\mathbb{C}$, $\mathbb{N}_{a+1}\triangleq\{a+1,a+2,a+3,\cdots\}$, $a\in\mathbb{R}$. More exactly, such a transform should be called the sampling free nabla Laplace transform, since the sampling period is assumed to be 1 \cite{Wei:2019FDTA}.

The region of convergence for $F(s)={{\mathscr N}_a}\left\{ {f\left( k \right)} \right\} $ is defined as the set of points in the complex plane for which the infinite series converges \cite{Goodrich:2015Book}
\begin{eqnarray}\label{Eq2}
{\textstyle {\rm{ROC}} \triangleq \left\{ {s:\big| {\sum\nolimits_{k = 1}^{ + \infty } {{{\left( {1 - s} \right)}^{k - 1}}f\left( {k + a} \right)} } \big| <  + \infty } \right\} .}
\end{eqnarray}

The inverse nabla Laplace transform can be written as \cite{Wei:2019AJC}
\begin{equation}\label{Eq3}
{\textstyle {\mathscr N}_a^{-1}\left\{ {F\left( s \right)} \right\}  \triangleq  \frac{1}{{2\pi {\rm{j}}}}\oint_c {F\left( s \right){{(1 - s)}^{ - k + a}}{\rm{d}}s} ,}
\end{equation}
where $k \in {\mathbb{N}_{a + 1}}$, $c$ is a closed curve rotating around the point $(1,{\rm j}0)$ clockwise and it also locates in the convergent region of $F(s)$ (defined as (\ref{Eq2})).

In this paper, the objective is to develop some methods to obtain the sequence $f(k)$ with $k\in\mathbb{N}_{a+1}$ from its nabla Laplace transform $F(s)={\mathscr N}_a\left\{ {f\left( k \right)} \right\} $ and the corresponding convergent region ${\rm ROC}$.

\section{MAIN RESULTS}\label{Section 3}
This section extends the conventional residual calculation method and partial fraction method to nabla Laplace domain. After carefully providing these methods, the advantage and disadvantage are then evaluated.

\subsection{Residual calculation method}
Based on the residue theorem \cite{Stein:2010Book}, the contour integral in equation (\ref{Eq3}) can be calculated accordingly.

When the finite poles of ${F\left( s \right){{\left( {1 - s} \right)}^{ - k + a}}}$ inside of the circle $c$ is $\{s_m\}$, then the sequence $f(k)$ with $k\in\mathbb{N}_{a+1}$ can be computed by the opposite number of the sum of all residues at the poles $\{s_m\}$, i.e.,
\begin{equation}\label{Eq4}
{\textstyle f\left( k \right) =  - \sum\limits_m {{\rm{Res}}\big[ {F\left( s \right){{\left( {1 - s} \right)}^{ - k + a}},{s_m}} \big]} .}
\end{equation}

If $s_m$ is $N$ times pole with $N\in\mathbb{Z}_+$, then the residue at $s_m$ is equal to
\begin{eqnarray}\label{Eq5}
{\textstyle
\begin{array}{l}
{\rm{Res}}\big[ {F\left( s \right){{\left( {1 - s} \right)}^{ - k + a}},{s_m}} \big] \\
= \frac{1}{{\left( {N - 1} \right)!}}\mathop {\lim }\limits_{s \to {s_m}} {\frac{{{{\rm{d}}^{N - 1}}}}{{{\rm{d}}{s^{N - 1}}}}{{\left( {s - {s_m}} \right)}^N}F\left( s \right){{\left( {1 - s} \right)}^{ - k + a}}}
\end{array}.}
\end{eqnarray}

If $s_m$ is single pole, then
\begin{eqnarray}\label{Eq6}
{\textstyle
\begin{array}{l}
{\rm{Res}}\big[ {F\left( s \right){{\left( {1 - s} \right)}^{ - k + a}},{s_m}} \big] \\
= \mathop {\lim }\limits_{s \to {s_m}} \left( {s - {s_m}} \right)F\left( s \right){\left( {1 - s} \right)^{ - k + a}}
\end{array}.}
\end{eqnarray}

When the finite poles of ${F\left( s \right){{\left( {1 - s} \right)}^{ - k + a}}}$ outside of the circle $c$ is $\{s_n\}$, then the sequence $f(k)$ with $k\in\mathbb{N}_{a+1}$ can be computed by
\begin{equation}\label{Eq7}
{\textstyle f\left( k \right) = \sum\limits_n {{\rm{Res}}\big[ {F\left( s \right){{\left( {1 - s} \right)}^{ - k + a}},{s_n}} \big]},}
\end{equation}
where the residue at the poles $\{s_n\}$ can be calculated like the previous discussion.

Note that the formulae in (\ref{Eq4}) and (\ref{Eq7}) are different from those in inverse $Z$ transform, since the clockwise curve $c$ is inverted. As a result, the sign should be specially handled when using the Residual theorem. It is worth emphasizing that for a finite value sequence $f(k)$, $s=1$ cannot be the pole of $F(s)$, since $f\left( {a + 1} \right) = \mathop {\lim }\limits_{s \to 1} F\left( s \right)\neq \infty$ can be obtained from the initial value theorem \cite{Wei:2019FDTA}.

\subsection{Partial fraction expansion method}

If the considered function $F(s)$ can be expressed as the following fractions
\begin{equation}\label{Eq8}
{\textstyle
F\left( s \right) = \sum\nolimits_{i = 1}^n {\frac{{{r_i}}}{{s - {s_i}}}} ,}
\end{equation}
where the coefficient ${r_i} = \mathop {\lim }\limits_{s \to {s_i}} \left( {s - {s_i}} \right)F\left( s \right)$, then the corresponding sequence satisfies
\begin{equation}\label{Eq9}
{\textstyle
f\left( k \right) = \sum\nolimits_{i = 1}^n {\frac{{{r_i}}}{{{{\left( {1 - {s_i}} \right)}^{k - a}}}}} ,k\in\mathbb{N}_{a+1}.}
\end{equation}

In equation (\ref{Eq7}), $F(s)$ has single pole $s_i$, $i=1,2,\cdots,n$. If there exists a multiple pole $\lambda$, i.e.,
\begin{equation}\label{Eq10}
{\textstyle
F\left( s \right) = \sum\nolimits_{i = 1}^{n - N} {\frac{{{r_i}}}{{s - {s_i}}}}  + \sum\nolimits_{i = 1}^N {\frac{{{q_i}}}{{{{\left( {s - \lambda } \right)}^i}}}},}
\end{equation}
where $q_i$, $i=1,2,\cdots,N$ can be calculated via ${q_i} = \mathop {\lim }\limits_{s \to \lambda } \frac{1}{{\left( {N - i} \right)!}}\frac{{{{\rm{d}}^{N - i}}}}{{{\rm{d}}{s^{N - i}}}}{\left( {s - \lambda } \right)^N}F\left( s \right)$, then the sequence $f(k)$, $k\in\mathbb{N}_{a+1}$ follows
\begin{eqnarray}\label{Eq11}
{\textstyle
f\left( k \right) = \sum\nolimits_{i = 1}^{n - N} {\frac{{{r_i}}}{{{{\left( {1 - {s_i}} \right)}^{k - a}}}}}  + \sum\nolimits_{i = 1}^N {\frac{{{q_i}{{\left( {k - a} \right)}^{i - 1}}}}{{\left( {i - 1} \right)!{{\left( {1 - \lambda } \right)}^{k - a + i - 1}}}}} .}
\end{eqnarray}

Likewise, the function $F(s)$ can be written as the sum of some familiar items, then the desired sequence can be obtained. For example, if the following equation holds
\begin{equation}\label{Eq12}
{\textstyle
F\left( s \right) = \sum\nolimits_{i = 1}^n {\frac{{{r_i}{s^{{\alpha _i} - {\beta _i}}}}}{{{s^{{\alpha _i}}} - {s_i}}}} ,}
\end{equation}
then one has
\begin{equation}\label{Eq13}
{\textstyle
f\left( k \right) = \sum\nolimits_{i = 1}^n {{r_i}{{\mathcal F}_{{\alpha _i},{\beta _i}}}\left( {{s_i},k,a} \right)} ,}
\end{equation}
where $\alpha_i\in\mathbb{R}$, $\beta_i\in\mathbb{R}_+$ and $k\in\mathbb{N}_{a+1}.$

This method can also be seen as the look-up table method. For convenience, with the help of some fundamental properties of nabla Laplace transform in \cite{Wei:2019FDTA}, 16 commonly used sequences and their nabla Laplace transforms are provided in Table \ref{Table 1}.

\begin{table*}
  \centering
  \caption{Nabla Laplace Transform Pairs.}
  \label{Table 1}
  \begin{tabular}{c l l l}
    \hhline
    number&$f(k)$, $k\in\mathbb{N}_{a+1}$, $a\in\mathbb{R}$        & $F(s)={\mathscr N}_a\left\{ {f\left( k \right)} \right\} $ & {\rm ROC}\\ \hline
    1&${\delta \left( {k-a - 1} \right)}$   & 1 & $s\in\mathbb{C}$\\ %\hline
    2&${u\left( {k-a - 1} \right)}$    &$\frac{{\rm{1}}}{s}$&$\left| {1 - s} \right| < 1$\\ %\hline
    3&${{k - a} }$&$\frac{{1}}{{{s^2}}}$&$\left| {1 - s} \right| < 1$\\ %\hline
    4&${{\gamma ^{k - a - 1}}}$&$\frac{1}{{1 - \gamma  + \gamma s}}$&$\left| {1 - s} \right|\left| \gamma  \right| < 1$, $\gamma\neq0$\\ %\hline
    5&$\frac{{{\left( {k - a} \right)}^{ \overline \alpha }}}{\Gamma \left( {\alpha  + 1} \right)}$&$\frac{1}{{{s^{\alpha  + 1}}}}$&$\left| {1 - s} \right| < 1,\alpha  \in \mathbb{C},\alpha  \notin {Z_- }$\\ %\hline
    6&${{\gamma ^{k - a - 1}}}\frac{{\left( {k - a} \right)}^{ \overline \alpha }}{\Gamma \left( {\alpha  + 1} \right)}$&$\frac{1}{{{{\left( {1 - \gamma  + \gamma s} \right)}^{\alpha  + 1}}}}$&$\left| {1 - s} \right|\left| \gamma  \right| < 1$, $\gamma\neq0$\\ %\hline
    7&${\frac{1}{{{{\left( {1 - \lambda } \right)}^{k - a}}}}}$&$\frac{1}{{s - \lambda }}$&$\left| {1 - s} \right| < \left| {1 - \lambda } \right|,\lambda  \ne 1$\\ %\hline
    8&$\frac{{{{\left( {k - a} \right)}^{ \overline {N-1} }}}}{{\left( {N - 1} \right)!{{\left( {1 - \lambda } \right)}^{k - a + N - 1}}}}$&$\frac{1}{{{{\left( {s - \lambda } \right)}^N}}}$&$\left| {1 - s} \right| < \min \left\{ {\left| {1 - \lambda } \right|,1} \right\},\lambda  \ne 1$, $N\in\mathbb{Z}_+$\\ %\hline
    9&${{{\mathcal F}_{\alpha ,\beta }}\left( {\lambda ,k,a} \right)}$&$\frac{{{s^{\alpha  - \beta }}}}{{{s^\alpha } - \lambda }}$&$\left| {1 - s} \right| < 1,\left| \lambda  \right| < {\left| s \right|^\alpha }$, $\alpha,\beta  \in {\mathbb{R}_ + }$\\ %\hline
    10&${\left( {k - a - 1} \right){{\mathcal F}_{\alpha ,\alpha }}\left( {\lambda ,k,a} \right)}$&$\frac{{\alpha {s^{\alpha  - 1}}\left( {1 - s} \right)}}{{{{\left( {{s^\alpha } - \lambda } \right)}^2}}}$&$\left| {1 - s} \right| < 1,\left| \lambda  \right| < {\left| s \right|^\alpha },\alpha  \in \mathbb{R}_+$\\ %\hline
    11&${{{\rm{e}}^{ - \lambda (k - a - 1)}}}$&$\frac{1}{{1 - {{\rm{e}}^{ - \lambda }}\left( {1 - s} \right)}}$&$\left| {1 - s} \right| < {{\rm{e}}^\lambda }$\\ %\hline
    12&${{\gamma ^{k - a - 1}}{{\rm{e}}^{ - \lambda (k - a - 1)}}}$&$\frac{1}{{1 - \gamma {{\rm{e}}^{ - \lambda }}\left( {1 - s} \right)}}$&$\left| {1 - s} \right|\left| \gamma  \right| < {{\rm{e}}^\lambda },\gamma  \ne 0$\\ %\hline
    13&${\sin \left( {\omega (k - a - 1)} \right)}$&$\frac{{\sin \left( \omega  \right)\left( {1 - s} \right)}}{{1 - 2\cos \left( \omega  \right)\left( {1 - s} \right) + {{\left( {1 - s} \right)}^2}}}$&$\left| {1 - s} \right| < 1$\\ %\hline
    14&${\cos \left( {\omega (k - a - 1)} \right)}$&$\frac{{1 - \cos \left( \omega  \right)\left( {1 - s} \right)}}{{1 - 2\cos \left( \omega  \right)\left( {1 - s} \right) + {{\left( {1 - s} \right)}^2}}}$&$\left| {1 - s} \right| < 1$\\ %\hline
    15&${\sinh \left( {\omega (k - a - 1)} \right)}$&$\frac{{\sinh \left( \omega  \right)\left( {1 - s} \right)}}{{1 - 2\cosh \left( \omega  \right)\left( {1 - s} \right) + {{\left( {1 - s} \right)}^2}}}$&$\left| {1 - s} \right| < \min \left\{ {{{\rm{e}}^\omega },{{\rm{e}}^{ - \omega }}} \right\}$\\ %\hline
    16&${\cosh \left( {\omega (k - a - 1)} \right)}$&$\frac{{1 - \cosh \left( \omega  \right)\left( {1 - s} \right)}}{{1 - 2\cosh \left( \omega  \right)\left( {1 - s} \right) + {{\left( {1 - s} \right)}^2}}}$&$\left| {1 - s} \right| < \min \left\{ {{{\rm{e}}^\omega },{{\rm{e}}^{ - \omega }}} \right\}$\vspace{3pt}\\ %\hline
    \hhline
  \end{tabular}
\end{table*}

A brief statement should be made in advance to facilitate understanding of Table \ref{Table 1}. $\delta \left( {n} \right) \triangleq \left\{ \begin{array}{l}
1,n = 0\\
0,n \ne 0
\end{array} \right.$ is the discrete-time unit impulse function.
$u\left( {n} \right) \triangleq \left\{ \begin{array}{l}
1,n \ge 0\\
0,n < 0
\end{array} \right.$ is the discrete-time unit step function.
${(k-a)}^{\overline{\alpha}}\triangleq\frac{\Gamma(k-a+\alpha)}{\Gamma(k-a)}$ is the rising function. ${{\mathcal F}_{\alpha ,\beta }}\left( {\lambda ,k,a} \right) \triangleq \sum\nolimits_{i = 0}^{ + \infty } {\frac{{{\lambda ^i}}}{{\Gamma \left( {i\alpha  + \beta } \right)}}{{\left( {k - a} \right)}^{\overline {i\alpha  + \beta  - 1} }}} $ is the discrete-time Mittag--Leffler function. ${\rm sinh}(\cdot)$ is the hyperbolic sine function and ${\rm cosh}(\cdot)$ is the hyperbolic cosine function. Besides, it can be calculated that ${{\mathscr N}_a}\left\{ {u\left( {k-a} \right)} \right\} = \frac{1}{s}={{\mathscr N}_a}\left\{ {u\left( {k-a-1} \right)} \right\} $, since we can only obtain the value of $f(k)$ with $k\in\mathbb{N}_{a+1}$ from $F(s)$ and ${u\left( {k-a} \right)}$ is exactly equal to ${u\left( {k-a-1} \right)}$ for any $k\in\mathbb{N}_{a+1}$.

For the purpose of comparison, more results are provided in Table \ref{Table 2} and Table \ref{Table 3}. More especially, Table \ref{Table 2} gives the generalized Laplace transform pairs. The transform and its inverse transform are defined as \cite{Wei:2018ArXivb}
\begin{equation}\label{Eq14}
{\textstyle {{\mathscr L}_a}\left\{ {f\left( t \right)} \right\} \triangleq \int_a^{ + \infty } {{{\rm{e}}^{ - s\left( {t - a} \right)}}f\left( t \right){\rm{d}}t} ,}
\end{equation}
\begin{equation}\label{Eq15}
{\textstyle {{\mathscr L}_a^{-1}}\left\{ {F\left( s \right)} \right\}  \triangleq \frac{1}{{2\pi {\rm{i}}}}\int_{\beta  - {\rm{i}}\infty }^{\beta  + {\rm{i}}\infty } { {{\rm{e}}^{s\left( {t - a} \right)}}{F\left( s \right)}{\rm{d}}s},}
\end{equation}
where $\beta$ is a real number so that the contour path of integration is in the convergence region of $F(s)$. When the sampling time $h$ is introduced in the nabla Laplace transform, the relationship between ${{\mathscr N}_a}\left\{ \cdot \right\}$ and ${{\mathscr L}_a}\left\{ \cdot \right\}$ can be derived (see \cite{Cheng:2011Book,Ortigueira:2016CNSNS}). The continuous-time Mittag--Leffler function is defined as ${{\mathcal E}_{\alpha ,\beta }}\left( {\lambda ,t,a} \right) \triangleq \sum\nolimits_{i = 0}^{ + \infty } {\frac{{{\lambda ^i}}}{{\Gamma \left( {i\alpha  + \beta } \right)}}{{\left( {t - a} \right)}^{i\alpha }}} $.

Table \ref{Table 3} gives the generalized Z-transform pairs. The transform and its inverse transform are defined as
\begin{equation}\label{Eq16}
{\textstyle {\mathscr Z}_a \left\{ {f\left( {k} \right)} \right\} \triangleq \sum\nolimits_{k = 0}^{ + \infty } {z^{-k}f\left( {k+a} \right) } , }
\end{equation}
\begin{equation}\label{Eq17}
{\textstyle {\mathscr Z}_a^{ - 1}\left\{ {F\left( z \right)} \right\} \triangleq \frac{1}{{2\pi {\rm{i}}}}\oint_c {{z^{k - a - 1}}F\left( z \right){\rm{d}}z} ,}
\end{equation}
where $c$ is a closed curve rotating around the point $(0,{\rm j}0)$ anticlockwise located in the convergent domain of $F(z)$ (see \cite{Wei:2019AJC}). By defining $g(k)=f(k+1)$, $k\in\mathbb{N}_{a}$ and $z^{-1}=1-s$, then ${\mathscr Z}_a \left\{ {g\left( {k} \right)} \right\} ={\mathscr N}_a \left\{ {f\left( {k} \right)} \right\} $.

\begin{table*}
  \centering
  \caption{Generalized Laplace Transform Pairs.}
  \label{Table 2}
  \begin{tabular}{c l l l}
    \hhline
    number&$f(t)$, $t\ge a$, $a\in\mathbb{R}$        & $F(s)={\mathscr L}_a\left\{ {f\left( t \right)} \right\} $ & {\rm ROC}\\ \hline
    1&${\delta \left( {t - a} \right)}$   & 1 & $s\in\mathbb{C}$\\ %\hline
    2&${u\left( {t-a} \right)}$    &$\frac{{\rm{1}}}{s}$&${\rm Re} \left\{ s \right\} > 0$\\ %\hline
    3&${{t - a} }$&$\frac{{1}}{{{s^2}}}$&${\rm Re} \left\{ s \right\} > 0$\\ %\hline
    4&${{\gamma ^{t - a}}}$&$\frac{1}{{ s-\ln\gamma}}$&${\rm Re} \left\{ s \right\} > \ln \gamma$, $\gamma  > 0$\\ %\hline
    5&$\frac{{{\left( {t - a} \right)}^{\alpha }}}{\Gamma \left( {\alpha  + 1} \right)}$&$\frac{1}{{{s^{\alpha  + 1}}}}$&${\mathop{\rm Re}\nolimits} \left\{ s \right\} > 0$, ${\rm Re} \left\{ \alpha  \right\} >  - 1$\\ %\hline
    6&${{\gamma ^{t - a }}}\frac{{\left( {t - a} \right)}^{\alpha }}{\Gamma \left( {\alpha  + 1} \right)}$&$\frac{1}{{{{\left( {s - \ln\gamma} \right)}^{\alpha  + 1}}}}$&${\rm Re} \left\{ s \right\} > \ln \gamma $, $\gamma  > 0$, ${\rm Re}\left\{ \alpha  \right\} >  - 1$\\ %\hline
    7&${{{\rm{e}}^{ \lambda (t - a)}}}$&$\frac{1}{{s - \lambda }}$&${\rm Re} \left\{ s \right\} > \lambda $\\ %\hline
    8&$\frac{{{{\left( {t- a} \right)}^{N - 1}}}}{{\left( {N - 1} \right)!}}{{\rm{e}}^{ \lambda (t - a)}}$&$\frac{1}{{{{\left( {s - \lambda } \right)}^N}}}$&${\rm Re} \left\{ s \right\} > \lambda$, $N\in\mathbb{Z}_+$\\ %\hline
    9&${\left( {t - a } \right)^{\beta-1}{{\mathcal E}_{\alpha ,\beta }}\left( {\lambda ,t,a} \right)}$&$\frac{{{s^{\alpha  - \beta }}}}{{{s^\alpha } - \lambda }}$&${\rm Re} \left\{ s \right\} > 0$, $\left| \lambda  \right| < {\left| s \right|^\alpha }$, $\alpha ,\beta  \in {\mathbb{R}_ + }$\\ %\hline
    10&${\left( {t - a } \right)^\alpha{{\mathcal E}_{\alpha ,\alpha }}\left( {\lambda ,t,a} \right)}$&$\frac{{\alpha {s^{\alpha  - 1}}}}{{{{\left( {{s^\alpha } - \lambda } \right)}^2}}}$&${\rm Re} \left\{ s \right\} > 0$, $\left| \lambda  \right| < {\left| s \right|^\alpha }$, $\alpha \in {\mathbb{R}_ + }$\\ %\hline
    11&${{{\rm{e}}^{ - \lambda (t - a )}}}$&$\frac{1}{s+\lambda}$&${\rm Re} \left\{ s \right\} > -\lambda $\\ %\hline
    12&${{\gamma ^{t - a }}{{\rm{e}}^{ - \lambda (t - a)}}}$&$\frac{1}{s+\lambda-\ln \gamma}$&${\rm Re} \left\{ s \right\} > \ln \gamma  - \lambda $, $\gamma  > 0$\\ %\hline
    13&${\sin \left( {\omega (t - a )} \right)}$&$\frac{\omega}{s^2+\omega^2}$&${\rm Re} \left\{ s \right\} > 0$\\ %\hline
    14&${\cos \left( {\omega (t - a)} \right)}$&$\frac{s}{s^2+\omega^2}$&${\rm Re} \left\{ s \right\} > 0$\\ %\hline
    15&${\sinh \left( {\omega (t - a)} \right)}$&$\frac{\omega}{s^2-\omega^2}$&${\rm Re} \left\{ s \right\} > \left| \omega  \right|$\\ %\hline
    16&${\cosh \left( {\omega (t - a)} \right)}$&$\frac{s}{s^2-\omega^2}$&${\rm Re} \left\{ s \right\} > \left| \omega  \right|$\vspace{3pt}\\ %\hline
    \hhline
  \end{tabular}
\end{table*}

\begin{table*}
  \centering
  \caption{Generalized Z-Transform Pairs.}
  \label{Table 3}
  \begin{tabular}{c l l l}
    \hhline
    number&$f(k)$, $k\in\mathbb{N}_{a+1}$, $a\in\mathbb{R}$        & $F(z)={\mathscr Z}_a\left\{ {f\left( k \right)} \right\} $ & {\rm ROC}\\ \hline
    1&${\delta \left( {k-a } \right)}$   & 1 & $s\in\mathbb{C}$\\ %\hline
    2&${u\left( {k-a } \right)}$    &$\frac{{\rm{1}}}{1-z^{-1}}$&$| z^{-1}| < 1$\\ %\hline
    3&${{k - a+1} }$&$\frac{1}{{{(1-z^{-1})^2}}}$&$| z^{-1}| < 1$\\ %\hline
    4&${{\gamma ^{k - a}}}$&$\frac{1}{{1 - \gamma z^{-1}}}$&$| z^{-1}|\left| \gamma  \right| < 1$, $\gamma\neq0$\\ %\hline
    5&$\frac{{{\left( {k - a+1} \right)}^{\bar \alpha }}}{\Gamma \left( {\alpha  + 1} \right)}$&$\frac{1}{{{(1-z^{-1})^{\alpha  + 1}}}}$&$| z^{-1}| < 1,\alpha  \in \mathbb{C},\alpha  \notin {Z_- }$\\ %\hline
    6&${{\gamma ^{k - a }}}\frac{{\left( {k - a+1} \right)}^{\bar \alpha }}{\Gamma \left( {\alpha  + 1} \right)}$&$\frac{1}{{{{( {1 - \gamma z^{-1}} )}^{\alpha  + 1}}}}$&$| z^{-1}|\left| \gamma  \right| < 1$, $\gamma\neq0$\\ %\hline
    7&${\frac{1}{{{{\left( {1 - \lambda } \right)}^{k - a+1}}}}}$&$\frac{1}{{1-z^{-1} - \lambda }}$&$| z^{-1}| < \left| {1 - \lambda } \right|,\lambda  \ne 1$\\ %\hline
    8&$\frac{{{{\left( {k - a+1} \right)}^{{\overline {^{N - 1}} }}}}}{{\left( {N - 1} \right)!{{\left( {1 - \lambda } \right)}^{k - a + N}}}}$&$\frac{1}{{{{( {1-z^{-1} - \lambda } )}^N}}}$&$| z^{-1}| < \min \left\{ {\left| {1 - \lambda } \right|,1} \right\},\lambda  \ne 1$, $N\in\mathbb{Z}_+$\\ %\hline
    9&${{{\mathcal F}_{\alpha ,\beta }}\left( {\lambda ,k+1,a} \right)}$&$\frac{{{(1-z^{-1})^{\alpha  - \beta }}}}{{{(1-z^{-1})^\alpha } - \lambda }}$&$| z^{-1}| < 1,\left| \lambda  \right| < {| 1-z^{-1}|^\alpha },\alpha ,\beta  \in {\mathbb{R}_ + }$\\ %\hline
    10&${\left( {k - a } \right){{\mathcal F}_{\alpha ,\alpha }}\left( {\lambda ,k+1,a} \right)}$&$\frac{{\alpha {(1-z^{-1})^{\alpha  - 1}}z^{-1}}}{{{{[ {{(1-z^{-1})^\alpha } - \lambda } ]}^2}}}$&$| z^{-1}| < 1,\left| \lambda  \right| < {| 1-z^{-1}|^\alpha },\alpha  \in \mathbb{R}_+$\\ %\hline
    11&${{{\rm{e}}^{ - \lambda (k - a)}}}$&$\frac{1}{{1 - {{\rm{e}}^{ - \lambda }}z^{-1}}}$&$| z^{-1}| < {{\rm{e}}^\lambda }$\\ %\hline
    12&${{\gamma ^{k - a }}{{\rm{e}}^{ - \lambda (k - a)}}}$&$\frac{1}{{1 - \gamma {{\rm{e}}^{ - \lambda }}z^{-1}}}$&$| z^{-1}|\left| \gamma  \right| < {{\rm{e}}^\lambda },\gamma  \ne 0$\\ %\hline
    13&${\sin \left( {\omega (k - a)} \right)}$&$\frac{{\sin \left( \omega  \right)z^{-1}}}{{1 - 2\cos \left( \omega  \right)z^{-1} + z^{-2}}}$&$| z^{-1}| < 1$\\ %\hline
    14&${\cos \left( {\omega (k - a)} \right)}$&$\frac{{1 - \cos \left( \omega  \right)z^{-1}}}{{1 - 2\cos \left( \omega  \right)z^{-1} + z^{-2}}}$&$| z^{-1}| < 1$\\ %\hline
    15&${\sinh \left( {\omega (k - a )} \right)}$&$\frac{{\sinh \left( \omega  \right)z^{-1}}}{{1 - 2\cosh \left( \omega  \right)z^{-1} + z^{-2}}}$&$| z^{-1}| < \min \left\{ {{{\rm{e}}^\omega },{{\rm{e}}^{ - \omega }}} \right\}$\\ %\hline
    16&${\cosh \left( {\omega (k - a)} \right)}$&$\frac{{1 - \cosh \left( \omega  \right)z^{-1}}}{{1 - 2\cosh \left( \omega  \right)z^{-1} + z^{-2}}}$&$| z^{-1}| < \min \left\{ {{{\rm{e}}^\omega },{{\rm{e}}^{ - \omega }}} \right\}$\vspace{3pt}\\ %\hline
    \hhline
  \end{tabular}
\end{table*}

\subsection{Suitability and limitations}
Both of the two developed methods can get the exact inverse nabla Laplace transform of some given $F(s)$ with special ROC. Nevertheless, every coin has two sides. They are generally used for $F(s)$ with rational polynomial expression and therefore their applicability of these methods should be pointed out clearly.

Nowadays, fractional calculus has played a critical role in many theoretical and practical scenarios. For this, when fractional sum or fractional difference of a sequence $f(k)$ is considered, the fractional order polynomial of the nabla Laplacian operator `$s$' will emerge in the expression of $F(s)$, such as $\frac{1}{s^\alpha-\lambda}$, $\frac{1}{(s-\lambda)^\alpha}$ and $\frac{1}{(s^2+1)^\alpha}$.

For the residual calculation method, if $F(s)=\frac{1}{s^\alpha-\lambda}$ with $\alpha\in(0,1)$, then there are infinite poles outside for the closed curve $c$ especially for the irrational case of $\alpha$ and the multiple pole $s=1$ has $k-a$ times. Therefore, it is difficult to compute $f(k)$ with the poles inside or outside $c$. Likewise, if $F(s)=\frac{1}{(s-\lambda)^\alpha}$ or $F(s)=\frac{1}{(s^2+1)^\alpha}$ with $\alpha\in(0,1)$ is considered, the problem follows immediately, i.e., $\lambda$ is $\alpha$ times pole of $F(s)$. As a result, both (\ref{Eq5}) and (\ref{Eq6}) cannot work for this case.

By adopting the basic properties of nabla Laplace transform, such as linearity, time advance, time delay, right shifting, left shifting, scaling in the frequency domain, differentiation in the frequency domain, integration in the frequency domain, accumulation, convolution and multiplication, etc, some transform pairs $f(k)\leftrightarrow F(s)$ can be found. Then, combining the obtained transform pairs like those in Table \ref{Table 1} with the partial fraction method, some inverse transform of fractional order polynomial.

Moreover, some special functions like $\frac{1}{{\rm e}^{s}-\lambda}$, $\frac{1}{{\rm log}(s^2+\lambda)}$, $\frac{\Gamma(s)}{\Gamma(s+\lambda)}$, $\frac{1}{{\rm tan}(1/s)}$, $\frac{1}{{\rm tanh}(s)}$ and $\frac{1}{{\rm sinh}(\sqrt{s})\sqrt{s}}$ will also appear under particular circumstances. The proposed two methods In other words, it is an opportunity and challenge to handle with such complicated irrational $F(s)$. Maybe, analytical or numerical techniques for the Laplace transform inversion \cite{Cohen:2007Book} could give us a lot of inspiration.

\section{EXAMPLES STUDY}\label{Section 4}
In this section, two illustrative examples are presented to further evaluate the theoretical approaches.

{\bf Example 1.} $F\left( s \right) = \frac{9}{{{{\left( {s + 1} \right)}^2}\left( {s - 2} \right)}}$, $\left| {1 - s} \right| < 1$.

\begin{enumerate}[(i)]
  \item Compute $f(k)$, $k\in\mathbb{N}_{a+1}$ with the residual calculation method. The poles of $F\left( s \right){\left( {1 - s} \right)^{ - k+a}}$ can be found as ${s_1} =  2,{s_2} = -1,{s_3} = 1$, where $s_1$ is a single pole, $s_2$ is a 2 times pole, and $s_3$ is a $k-a$ times pole. Selecting a closed curve $c$ in the convergent region $\left| {1 - s} \right| < 1$, then $s_1$, $s_2$ are outside the curve and $s_3$ is inside the curve.

      By applying formula (\ref{Eq4}), it follows
      \[
      \begin{array}{l}
        f\left( k \right)\\
         =  - {\rm{Res}}\big[ {F\left( s \right){{\left( {1 - s} \right)}^{- k +a }},1} \big]\\
         = \frac{{ - 1}}{{\left( {k - a - 1} \right)!}}\mathop {\lim }\limits_{s \to 1} \frac{{{{\rm{d}}^{k - a - 1}}}}{{{\rm{d}}{s^{k - a - 1}}}}\big[ {{{\left( {s - 1} \right)}^{k - a}}F\left( s \right){{\left( {1 - s} \right)}^{ - k+a}}} \big]\\
         = \frac{{{{\left( { - 1} \right)}^{k - a + 1}}}}{{\left( {k - a - 1} \right)!}}\mathop {\lim }\limits_{s \to 1} \frac{{{{\rm{d}}^{k - a - 1}}}}{{{\rm{d}}{s^{k - a - 1}}}}\big[ {\frac{1}{{s - 2}} - \frac{1}{{s + 1}} - \frac{3}{{{{\left( {s + 1} \right)}^2}}}} \big]\\
         = \mathop {\lim }\limits_{s \to 1} \big[ {\frac{1}{{{{\left( {s - 2} \right)}^{k - a}}}} - \frac{1}{{{{\left( {s + 1} \right)}^{k - a}}}} - \frac{{3\left( {k - a} \right)}}{{{{\left( {s + 1} \right)}^{k - a + 1}}}}} \big]\\
         = {\left( { - 1} \right)^{k - a}} - {2^{a - k}} + 3\left( {a - k} \right){2^{a - k - 1}}.
        \end{array}\]

    If we use formula (\ref{Eq7}), then
    \[\begin{array}{rl}
    f\left( k \right) =\hspace{-6pt}& {\rm{Res}}\big[ {F\left( s \right){{\left( {1 - s} \right)}^{ - k + a}},2} \big]\\
     \hspace{-6pt}&+ {\rm{Res}}\big[ {F\left( s \right){{\left( {1 - s} \right)}^{ - k + a}}, - 1} \big]\\
     =\hspace{-6pt}& \mathop {\lim }\limits_{s \to  2} \big[ {\left( {s - 2} \right)F\left( s \right){{\left( {1 - s} \right)}^{ - k + a}}} \big]\\
     \hspace{-6pt}&+ \mathop {\lim }\limits_{s \to - 1} \big[ {\left( {s + 1} \right)F\left( s \right){{\left( {1 - s} \right)}^{ - k + a}}} \big]\\
     =\hspace{-6pt}& {\left( { - 1} \right)^{k - a}} - {2^{a - k}} + 3\left( {a - k} \right){2^{a - k - 1}}.
    \end{array}\]

    \item Compute $f(k)$, $k\in\mathbb{N}_{a+1}$ with the partial fraction method. After simple mathematic deduction, one has
       \[{\textstyle F\left( s \right) = \frac{1}{{s - 2}} - \frac{1}{{s + 1}} - \frac{3}{{{{\left( {s + 1} \right)}^2}}}.}\]

        By the aid of Table \ref{Table 1}, one obtains
        \[\begin{array}{l}
        f\left( k \right)\\
         = {\mathscr N}_a^{ - 1}\big\{ {\frac{1}{{s - 2}}} \big\} - {\mathscr N}_a^{ - 1}\big\{ {\frac{1}{{s + 1}}} \big\} - {\mathscr N}_a^{ - 1}\big\{ {\frac{3}{{{{\left( {s + 1} \right)}^2}}}} \big\}\\
         = \frac{1}{{{{\left( {1 - 2} \right)}^{k - a}}}} - \frac{1}{{{{\left( {1 + 1} \right)}^{k - a}}}} - \frac{{3\left( {k - a} \right)}}{{{{\left( {1 + 1} \right)}^{k - a + 1}}}}\\
         = {\left( { - 1} \right)^{k - a}} - {2^{a - k}} + 3\left( {a - k} \right){2^{a - k - 1}}.
        \end{array}\]
\end{enumerate}

The first method seems more complex than the second one, since some arithmetic operations have been done when building  Table \ref{Table 1}. However, we have to admit that both the two established methods can solve the problem exactly.

{\bf Example 2.} $F\left( s \right) = \frac{{0.2{s^{0.2}} - 0.3}}{{{s^{1.2}} - 0.2{s^{0.7}} - 0.3{s^{0.5}} + 0.06}}$, $\left| {1 - s} \right| < 1$ and $\left|s \right|>0.3^{10/7}$.

\begin{enumerate}[(i)]
  \item Compute $f(k)$, $k\in\mathbb{N}_{a+1}$ with the residual calculation method. The poles of $F\left( s \right){\left( {1 - s} \right)^{ - k+a}}$ can be found as ${s_1} = 0.04,{s_2} = 0.3^{10/7}{{\rm{e}}^{{\rm{j}}20\pi i/7}},{s_3} = 1$, where $s_1$ is a single pole and $s_3$ is a $k-a$ times pole. Though $s_2$ is a single pole, different value of $i$ correspond to different $s_2$. In fact, the number of $s_2$ is infinite. Selecting a closed curve $c$ in the convergent region, then $s_1$, $s_2$ are outside the curve and $s_3$ is inside the curve. To avoid the infinite number of poles, we utilize the pole inside the curve $c$. Along this way, one has
      \[\begin{array}{l}
    f\left( k \right) \\
    =  - {\rm{Res}}\big[ {F\left( s \right){{\left( {1 - s} \right)}^{ - k + a}},1} \big]\\
    = \frac{{ - 1}}{{\left( {k - a - 1} \right)!}}\mathop {\lim }\limits_{s \to 1} \frac{{{{\rm{d}}^{k - a - 1}}}}{{{\rm{d}}{s^{k - a - 1}}}}\big[ {{{\left( {s - 1} \right)}^{k - a}}F\left( s \right){{\left( {1 - s} \right)}^{ - k + a}}} \big]\\
    = \frac{{{{\left( { - 1} \right)}^{k - a + 1}}}}{{\left( {k - a - 1} \right)!}}\mathop {\lim }\limits_{s \to 1} \frac{{{{\rm{d}}^{k - a - 1}}}}{{{\rm{d}}{s^{k - a - 1}}}}\big[\frac{{0.2{s^{0.2}} - 0.3}}{{{s^{1.2}} - 0.2{s^{0.7}} - 0.3{s^{0.5}} + 0.06}}\big]
    \end{array}\]

    To be honest, it is difficult to provide a brief expression of $\frac{{{{\rm{d}}^{k - a - 1}}}}{{{\rm{d}}{s^{k - a - 1}}}}\big[\frac{{0.2{s^{0.2}} - 0.3}}{{{s^{1.2}} - 0.2{s^{0.7}} - 0.3{s^{0.5}} + 0.06}}\big]$. In other words, the residual calculation method cannot solve the problem effectively.

    \item Compute $f(k)$, $k\in\mathbb{N}_{a+1}$ with the partial fraction method. The considered function $F(s)$ can be equivalently rewritten as
        \[{\textstyle F\left( s \right) = \frac{1}{{{s^{0.5}} - 0.2}} - \frac{{{s^{0.2}}}}{{{s^{0.7}} - 0.3}}}.\]

        Similarly, using Table \ref{Table 1} yields
        \[\begin{array}{rl}
        f\left( k \right) =\hspace{-6pt}& {\mathscr N}_a^{ - 1}\big\{ {\frac{1}{{{s^{0.5}} - 0.2}}} \big\} - {\mathscr N}_a^{ - 1}\big\{ {\frac{{{s^{0.2}}}}{{{s^{0.7}} - 0.3}}} \big\}\\
         =\hspace{-6pt}& {{\mathcal F}_{0.5,0.5}}\left( {0.2,k,a} \right) - {{\mathcal F}_{0.7,0.5}}\left( {0.3,k,a} \right),
        \end{array}\]
        which means that the partial fraction method can compute the sequence $f(k)$ from some special irrational $F(s)$. Generally speaking, it is essential to disassemble $F(s)$ into some essential elements and build a more detailed nabla Laplace pairs table include these elements.
\end{enumerate}

%\balance
\section{CONCLUSIONS}\label{Section 5}
In this paper, the conventional residual calculation method and partial fraction expansion method have been investigated for computing the nabla Laplace transform. It is the first time to give two available methods instead of calculating the contour integral directly. Although the exact solution can be achieved, the developed methods have limitations. More related methods are expected for the irrational case. It is hoped that this paper will be a useful tool for all those who use nabla Laplace transforms in their work.

\end{document}